\documentclass[reqno,fleqn]{amsart}%
\usepackage[english]{babel}
\usepackage{amsmath}
\usepackage{amssymb}
\usepackage{amsfonts}
\usepackage{graphicx}%
\theoremstyle{plain}
\newtheorem{defn}{Definition}[section]

\newtheorem{prop}[defn]{Proposition}

\newtheorem{thm}[defn]{Theorem}

\newtheorem{cor}[defn]{Corollary}

\newtheorem{rem}[defn]{Remark}

\numberwithin{equation}{section} \allowdisplaybreaks

\newcommand\R{\mathbb{R}}
\newcommand{\C}{\mathbb{C}}

\newcommand{\cal}{\mathcal }
\newcommand\g {\mathfrak{g}}
\newcommand\p {\mathfrak{p}}
\newcommand\ki {\mathfrak{k}}
\newcommand\so {\mathfrak{so}}
\newcommand\sln {\mathfrak{sl}}
\newcommand\n {\mathfrak{n}}

\begin{document}
\title{A geometric proof of the Karpelevich-Mostow's Theorem}
\author{Antonio J. Di Scala}
\address{A.J.D.: Dipartimento di Matematica, Politecnico di Torino, Corso Duca degli
Abruzzi 24, 10129 Torino, Italy }
\email{antonio.discala@polito.it}
\author{Carlos Olmos}
\address{C.O.: Fa.M.A.F., Universidad Nacional de Cordoba, Ciudad Universitaria,
5000 Cordoba, Argentina} \email{olmos@mate.uncor.edu}
\thanks{Research supported by Programa Raices, Subsidio Cesar Milstein, Republica Argentina.}
\date{\today}

\begin{abstract}
In this paper we give a geometric proof of the Karpelevich's theorem that asserts that a semisimple Lie subgroup
of isometries, of a symmetric space of non compact type, has a totally geodesic orbit. In fact, this is
equivalent to a well-known result of Mostow about existence of compatible Cartan decompositions.

\end{abstract}

\maketitle

\section{Introduction.}

In this paper we address the problem of giving a geometric proof of
the following theorem of Karpelevich.

\begin{thm}\label{Ka}(Karpelevich \cite{Ka})
Let $M$ be a Riemannian symmetric space of non positive curvature
without flat factor. Then any connected and semisimple subgroup $G
\subset \text{Iso}(M)$ has a totally geodesic orbit $G.p \subset M$.
\end{thm}

It is well-known that Karpelevich's theorem is equivalent to the
following algebraic theorem.

\begin{thm} \label{Mo}(Mostow \cite[Theorem 6]{Mo})
Let $\g'$ be a real semisimple Lie algebra of non compact type and
let $\g \subset \g'$ be a semisimple Lie subalgebra. Let $\g = \ki
\oplus \p$ be a Cartan decomposition for $\g$. Then there exists a
Cartan decomposition $\g' = \ki' \oplus \p'$ for $\g'$ such that
$\ki \subset \ki'$ and $\p \subset \p'$.
\end{thm}

The proof of the above theorems is very algebraic in nature and uses
delicate arguments related to automorphisms of semisimple Lie algebras
(see \cite{Onish} for a modern exposition of Karpelevich's results).\\

For the real hyperbolic spaces, i.e. when $\g'= \so(n,1)$, there are two geometric proofs of Karpelevich's
theorem \cite{olmos-discala01}, \cite{BZ}. The proof in \cite{olmos-discala01} is based on the study of minimal
orbits of isometry subgroups, i.e. orbits with zero mean curvature. The approach in \cite{BZ} is based on
technics from hyperbolic dynamics. It is interesting to note that both proofs are strongly based on the fact
that the boundary at infinity of real hyperbolic spaces has a
simple structure.\\

The only non-trivial algebraic tool that we will use is the
existence of a Cartan decomposition of a non compact semisimple Lie
algebra. But this can also be proved geometrically as
was explained by S.K. Donaldson in \cite{Do07}. \\

Here is a brief explanation of our proof of Theorem \ref{Ka}. We
first show that a simple subgroup $G \subset \text{Iso}(M)$  has a
minimal orbit $G.p \subset M$. Then, by using a standard totally
geodesic embedding $M \hookrightarrow \cal{P}$, where ${\cal P} =
SL(n,\R)/SO(n)$, we will show that $G.p$ is, actually,
a totally geodesic submanifold of $M$. \\

\section{Preliminaries.}

The equivalence between Theorems \ref{Ka} and \ref{Mo}, for the non
expert reader, is a consequence of the following Elie Cartan's
famous and remarkable theorem.

\begin{thm} (Elie Cartan) Let $M$ be a Riemannian symmetric space of
non positive curvature without flat factor. Then the Lie group
$\text{Iso}(M)$ is semisimple of non compact type. Conversely, if
$\g$ is a semisimple Lie algebra of non compact type then there
exist a Riemannian symmetric space $M$ of non positive curvature
without flat factor such that $\g$ is the Lie algebra of
$\text{Iso}(M)$.\\
\end{thm}

\begin{rem} \label{cd} \rm We remark, for the non expert reader, that the difficult
part of the proof of the above proposition is the second part.
Namely, the construction of the Cartan decomposition $\g = \ki
\oplus \p$, where $\ki$ is maximal compact subalgebra of $\g$ and
the Killing form $B$ of $\g$ is positive definite on $\p$. The
standard and well-known proof of the existence of a Cartan
decomposition is long and via the classification theory of complex
semisimple Lie algebras, i.e. the existence of a real compact form
(see e.g. \cite{Helgason78}). There is also a direct and geometric
proof of the existence of a Cartan decomposition \cite{Do07}.

On the other hand, when $\g = \text{Lie}(\text{Iso}(M))$, where $M$
is a Riemannian symmetric space of non positive curvature without
flat factor, a Cartan decomposition of $\g$ can be constructed
geometrically. Namely, $\g = \text{Lie}(\text{Iso}(M)) = \ki \oplus
\p$ where $\ki$ is the Lie algebra of the isotropy group $K_p
\subset \text{Iso}(M)$ and $\p : = \{ X \in
\text{Lie}(\text{Iso}(M)) :
(\nabla X)_p = 0 \}$.\\
\end{rem}

It is well-known that the Riemannian symmetric spaces ${\cal P} =
SL(n, \R)/SO(n)$ are the \it universal \rm Riemannian symmetric
space of non positive curvature. Namely, any Riemannian symmetric
space of non compact type $M = G/K$ can be totally geodesically
embedded in some ${\cal P}$ (up to rescaling the metric in the
irreducible De Rham factors). A proof of this fact follows from the
following well-known result (c.f. Theorem 1 in \cite{Do07}).

\begin{prop} \label{SL} Let $\g \subset \sln(n,\R)$ be a semisimple Lie
subalgebra and let $\g = \ki \oplus \p $ be a Cartan decomposition.
Then there exists a Cartan decomposition $\sln(n,\R) = {\cal A}
\oplus {\cal S}$ such that $\ki \subset {\cal A}$ and $\p \subset
{\cal S}$. Thus, if $G \subset SL(n)$ is semisimple, $G$ has a
totally geodesic orbit in ${\cal P} = SL(n)/SO(n)$. Indeed, any
Riemannian symmetric space of non positive curvature $M$, without
flat factor, can be totally geodesically embedded in some ${\cal P}
= SL(n)/SO(n)$.
\end{prop}

\it Proof. \rm Notice that any Cartan decomposition of $\sln(n,\R)$
is given by the anti-symmetric ${\cal A}$ and symmetric matrices
${\cal S}$ w.r.t. a positive definite inner product on $\R^n$. Since
$ g^* := \ki \oplus i \p$ is a compact Lie subalgebra of
$\sln(n,\C)$, there exists a positive definite Hermitian form $( \ |
\ )$ of $\C^n$ invariant by $ g^* $. By defining $\langle \, , \,
\rangle := Real ( \ | \ )$ it follows that $ \ki \subset {\cal A}$
and $ \p\subset {\cal S}$.
$\Box$\\

Here is another corollary of the existence of the totally geodesic
embedding $M \hookrightarrow {\cal P}$.

\begin{cor} \label{Bo} Let $M$ be a Riemannian symmetric space of non positive curvature
without flat factor and let $\partial M$ be its boundary at
infinity. Then a connected and semisimple Lie subgroup $G \subset
\text{Iso}(M)$ of non compact type has no fixed points in $\partial
M$.
\end{cor}

Finally, for the non expert reader, we include the following
proposition.

\begin{prop}  \label{euclidean} Let $M$ be a Riemannian symmetric space of
non positive curvature without flat factor. Let $S = \R^N \times M$
be a symmetric space of non positive curvature with flat factor
$\R^N$. If $G \subset \text{Iso}(S)$ is a connected non compact
simple Lie group then $G \subset \text{Iso}(M)$.
\end{prop}

\it Proof. \rm Let $\g$ be the Lie algebra of $G$. Then the
projection $\pi: \g \mapsto Lie(Iso(\R^N))$ is injective or trivial
i.e. $\pi \equiv 0$. If $\pi$ is injective then a further
composition with the projection to $so(N)$ gives that $\g$ must
carry a bi-invariant metric. So, $\g$ can not be simple and non
compact. $\Box$

\section{Minimal and totally geodesic orbits.}

We will need the following proposition (see Proposition 5.5. in
\cite{AlDi03} or Lemma 3.1. in \cite{olmos-discala01}).

\begin{prop} \label{u2} Let $M$ be a Riemannian symmetric space of non positive curvature
without flat factor and let $G \subset \text{Iso}(M)$ be a connected
group of isometries. Assume that $G$ has a totally geodesic orbit
$G.p$. Then any other minimal orbit $G.q $ is also a totally
geodesic submanifold of $M$.
\end{prop}

\it Proof. \rm Let $G.p$ be the totally geodesic orbit and let $G.q
\neq \{q\}$ be another orbit. Notice that $G.p$ is a closed embedded
submanifold of $M$. Let $\gamma$ be a geodesic in $M$ that minimizes
the distance between $q$ and $G.p$. Without loss of generality we
may assume that $\gamma(0)=p$ and $\gamma(1)=q$ (eventually by
changing the base point $p$ by another in the orbit). It is standard
to show that $\dot{\gamma}(0)$ is perpendicular to $T_p(G.p)$, or
equivalently $\langle X.p,\dot{\gamma}(0) \rangle = 0$ for all $X$
in the Lie algebra of $G$. Observe that this implies that $\langle
X.\gamma(t) , \dot{\gamma}(t) \rangle = 0$ for all $t$, because of
$\frac{d}{dt} \langle X.\gamma(t), \dot{\gamma}(t) \rangle
 = \langle \nabla_{\dot{\gamma}(t)}X.\gamma(t),
 \dot{\gamma}(t) \rangle = 0$, by the Killing equation.
 So, $\dot{\gamma}(t)$ is perpendicular to
$T_{\gamma(t)}(G.\gamma(t))$, for all $t$.
 Let $X$ be a Killing field in the Lie algebra of $G$
 such that $X.q \neq 0$  and let $\phi_s^X$ be
the one-parameter group of isometries generated by $X$. Define
$h:I\times\R \rightarrow M$ by $h_s(t):=\phi_s^X.\gamma(t)$. Note
that $X.h_s(t) = \frac{\partial h}{\partial s}$ and that, for a
fixed $s$, $h_s(t)$ is a geodesic.

Let $A_{\dot{\gamma}(t)}$ be the shape operator, in the direction of
$\dot{\gamma}(t)$ of the orbit $G.\gamma(t)$. Define $f(t) := -
\langle A_{\dot{\gamma}(t)}(X.\gamma(t)) ,X.\gamma(t) \rangle =
  \langle \frac{D}{\partial s}\frac{\partial h}{\partial t},
X.h_s(t) \rangle \mid_{s=0}$. Let us compute
$$\frac{d}{dt}f(t) = \langle \frac{D}{\partial t}
\frac{D}{\partial s}\frac{\partial h}{\partial t}, X.h_s(t) \rangle
\mid_{s=0} + \langle \frac{D}{\partial s} \frac{\partial h}{\partial
t}, \frac{D}{\partial t}X.h_s(t) \rangle \mid_{s=0} $$
$$  = \langle R(\frac{\partial h}{\partial t},
\frac{\partial h}{\partial s}) \frac{\partial h}{\partial t} ,
X.h_s(t) \rangle \mid_{s=0}
 + \langle \frac{D}{\partial t}\frac{\partial h}{\partial s}, \frac{D}{\partial t}X.h_s(t) \rangle \mid_{s=0} $$
$$ = \langle R(\dot{\gamma}(t),X.\gamma(t))\dot{\gamma}(t),
X.\gamma(t) \rangle + \| \nabla_{\dot{\gamma}(t)}(X.\gamma(t)) \| ^2
$$

Then $\frac{d}{dt}f(t) \geq 0$.
 Since $f(0)=0$, due to the fact that $G.p$
is totally geodesic, we obtain that $f(1)= - \langle
A_{\dot{\gamma}(1)}(X.q),X.q \rangle \geq 0$. Hence
$A_{\dot{\gamma}(1)}$ is negative semidefinite. Since $G.q$ is
minimal, $\text{trace}(A_{\dot{\gamma}(1)})=0$, we get that  $f(t)
\equiv 0$. Thus, $\langle
R(\dot{\gamma}(t),X.\gamma(t))\dot{\gamma}(t), X.\gamma(t) \rangle
\equiv 0$ and $\nabla_{\dot{\gamma}(t)}(X.\gamma(t)) \equiv 0$.
Notice that the tangent spaces $T_{\gamma(t)} G \cdot \gamma(t)$ are
parallel along $\gamma(t)$ in $M$. So the normal spaces
$\nu_{\gamma(t)} G \cdot \gamma(t)$ are also parallel along
$\gamma(t)$ in $M$. Since $M$ is a symmetric space of non positive
curvature the condition $\langle
R(\dot{\gamma}(t),X.\gamma(t))\dot{\gamma}(t), X.\gamma(t) \rangle
\equiv 0$ implies $R(\dot{\gamma}(t),X.\gamma(t))(\cdot) \equiv 0 $.
Let $\eta(t) \in \nu_{\gamma(t)} G.\gamma(t) $ be a parallel vector
along $\gamma(t)$ and let $X,Y$ two Killing vector fields in the Lie
algebra of $G$. Then $\frac{d}{dt} \langle \nabla_XY, \eta(t)
\rangle = \langle \nabla_{\dot{\gamma}(t)} \nabla_XY, \eta(t)
\rangle = \langle \nabla_X \nabla_{\dot{\gamma}(t)} Y, \eta(t)
\rangle + \langle R(\dot{\gamma}(t),X.\gamma(t))(Y), \eta(t) \rangle
\equiv 0$. Since $\langle \nabla_XY, \eta(0) \rangle = 0$ we get
that the $G$-orbits $G. \gamma(t)$ are totally geodesic submanifolds
of $M$. $\Box$ \\

Let $T_G := \{ p \in M : G.p \mbox{ is a totally geodesic
submanifold of } M \}$ be the union of totally geodesic $G$-orbits.

\begin{cor} \label{u3} In the assumption of the above proposition assume also
that $G$ is semisimple. Then $T_G$ is a totally geodesic submanifold
of $M$. Moreover $T_G$ is a Riemannian product $T_G = (G.p) \times
N$ where $N$ is a totally geodesic submanifold of $M$.
\end{cor}

\it Proof. \rm Let $K':= \text{Iso}(M)_p$ be the isotropy subgroup at $p \in M$ and let $\ki'$ its Lie algebra.
Let $\p' \subset \text{Lie}(\text{Iso}(M))$ such that $X \in \p'$ iff $(\nabla X)_p = 0$. Thus,
$\text{Lie}(\text{Iso}(M)) = \ki' \oplus \p'$ is a Cartan decomposition of $\text{Lie}(\text{Iso}(M))$. Let $\g
= \ki \oplus \p$ be a Cartan decomposition of the Lie algebra $\g = \text{Lie}(G)$. Since $G.p$ is totally
geodesic in $M$ we get that $\ki \subset \ki'$ and $\p \subset \p'$. Let $\n := \{ Y \in \p' : Y \perp \p \mbox{
and } [Y, \p] = 0 \}$ which is a Lie triple system of $\p'$. Moreover,  $\p \oplus \n$ is also a Lie triple
system of $\p'$. So, $exp_p (\p \oplus \n) = exp_p (\p) \times exp_p (\n) $ is a $G$-invariant totally geodesic
submanifold of $M$. Notice that (by construction) $exp_p (\p \oplus \n) \subset T_G$. Finally, the last part of
the proof of the above Proposition \ref{u2} shows that
$T_G \subset exp_p (\p \oplus \n) $.$\Box$\\

\section{Karpelevich's Theorem for $G$ a simple Lie group.}

Here is the first step to prove Theorem \ref{Ka}.

\begin{thm} Let $M$ be a Riemannian symmetric space of non positive curvature
without flat factor. Then any connected, simple and non compact Lie
subgroup $G \subset {Iso}(M)$ has a minimal orbit $G.p \subset M$.
\end{thm}

\it Proof. \rm Let $\g = \ki \oplus \p$ be a Cartan decomposition of the Lie algebra $ \g := Lie(G)$ and let $K
\subset G$ be the maximal compact subgroup associated to $\ki$. Let $\Sigma$ be the set of fixed points of $K$.
Notice that $\Sigma \neq \emptyset$ by Cartan's fixed point theorem.  Notice that $G$-orbits $G.x$ trough points
in $x \in \Sigma$ are irreducible Riemannian symmetric spaces associated to the symmetric pair $(G,K)$. So, if
$x,y \in \Sigma$, the metric induced on $G.x$ and $G.y$ differ from a constant multiple. Let $x_0 \in \Sigma$
and let $g_0$ be the Riemannian metric of $G.x_0 = G/K$ induced by the Riemannian metric $g = \langle,\rangle$
of $M$. So, if $y \in \Sigma$ the Riemannian metric $g_y$ on $G.y$ is given by $g = \lambda(y) \cdot g_0$.
Notice that if $X \in \p$ is unitary at $x_0$ (i.e. $g_0(X(x_0), X(x_0)) = 1 $) then $\lambda(y) = g(X(y),X(y))
=
 \|X(y)\|^2$. We claim that $\lambda(y)$ has a minimum in $\Sigma$.
Indeed, if $y_n \rightarrow \infty \in
\partial \Sigma \subset \partial M$ ($y_n \in \Sigma $) and $\lambda(y_n) \leq const$
then the monoparametric Lie group $\psi_t^X \subset G$ associated to
any unitary $X \in \p$ at $x_0 \in \Sigma$ must fix $\infty \in
\partial \Sigma \subset
\partial M$. Thus, since $X \in \p$ is arbitrary and $\p$ generate $\g$ we
get that $\infty \in
\partial \Sigma \subset \partial M$ is a fixed point of $G$. This
contradicts Proposition \ref{Bo}. So there exist $y_0 \in \Sigma$ such that $\lambda(y_0)$ is minimum. Notice
that the volume element $Vol_y$ of an orbit $G.{y}$ is given by ${\lambda}^{\frac{n}{2}} Vol_{x_0}$, where $n =
dim(G/K)$. Now a simple computation shows that the mean curvature vector of $G.{y_0}$ vanishes.
$\Box$\\

Now we are ready to prove Karpelevich's Theorem \ref{Ka} for $G$
a simple non compact Lie subgroup of $\text{Iso}(M)$.\\

\begin{thm} \label{KaS} Let $M$ be a Riemannian symmetric space of non positive curvature.
Then any connected, simple and non compact Lie subgroup $G \subset \text{Iso}(M)$ has a totally geodesic orbit
$G.p \subset M$.
\end{thm}

\it Proof. \rm According to Proposition \ref{euclidean} we can assume that $M$ has no flat factor. Let $i: (M,g)
\hookrightarrow (\cal{P}, h)$ be a totally geodesic embedding as in Proposition \ref{SL}. Notice that the
pull-back metric $i^*h$ can eventually differ (up to constants factors) from $g$ on each irreducible De Rham
factor of $M$. Anyway, totally geodesic submanifolds of $(M,g)$ and $(M,i^*h)$ are the same. Notice that $G$
also acts by isometries on $(M,i^*h)$. Indeed, $G$ can be regarded as a subgroup of $\text{Iso}(\cal P)$. Now
Proposition \ref{SL} implies that $G$ has a totally geodesic orbit $G.p$ in ${\cal P}$. The above proposition
shows that $G$ has a minimal orbit $G.y_0$ in $(M,i^*h)$. Since the embedding $M \hookrightarrow \cal{P}$ is
totally geodesic we get that the $G$-orbit $G.y_0$ is also a minimal submanifold of $\cal{P}$. Then Proposition
\ref{u2} implies that $G.y_0$ is a totally geodesic submanifold of $\cal{P}$. Thus, $G.y_0$ is a totally
geodesic submanifold of $(M,i^*h)$ and so $G.y_0$ is also a totally geodesic submanifold of $(M,g)$. $\Box$

\section{Karpelevich's Theorem.}

Let $G \subset \text{Iso}(M)$ be a semisimple, connected Lie group. Then the Lie algebra $\g = Lie(G) = \g_1
\oplus \g_2$ is a sum of a simple Lie algebra $\g_1$ and a semisimple Lie algebra $\g_2$. Due to Cartan's fixed
point theorem we can assume that each simple factor of $\g$ is non compact. We will make induction on the number
of simple factors of the semisimple Lie algebra $\g$. Let $G_1$ (resp. $G_2$) be the simple Lie group associated
to $\g_1$ (resp. the semisimple Lie subgroup associated to $\g_2$). Let $T_{G_1} \subset M$ be the union of the
totally geodesic orbits of the simple subgroup $G_1$ acting on $M$. Notice that Theorem \ref{KaS} implies that
$T_{G_1} \neq \emptyset$ and Proposition \ref{u3} implies that $T_{G_1} = (G_1 \cdot p) \times \, N$ is a
totally geodesic submanifold of $M$, where $G_1 \cdot p $ is a totally geodesic $G_1$-orbit. Notice that $G_2$
acts on $T_{G_1} = (G_1 \cdot p) \times \, N$, since $G_2$ commutes with $G_1$. Then $\g_2$ (or eventually a
quotient $\g_2/\sim$ of it) acts on $N$. Since $N$ is symmetric space of non positive curvature we get, by
induction, that the semisimple subgroup $G_2$ (or eventually a quotient $G_2/\sim$ of it) has a totally geodesic
orbit $ S \subset N $. Then $(G_1 \cdot p) \times S$ is a totally geodesic orbit of $G$ and this finish our
proof of Karpelevich's Theorem \ref{Ka}.

\end{document}